\documentclass[11pt]{amsart}

\usepackage{amssymb} 
\usepackage[latin1]{inputenc}
\usepackage[matrix,arrow,curve]{xy}
\usepackage[mathscr]{eucal}
\usepackage{epic,eepic}
\usepackage{bbm}

\makeindex
\newtheorem{theorem}{Theorem}

\newtheorem{prop}{Proposition}
\newtheorem{fact}{Fact}

\theoremstyle{definition}
\newtheorem{definition}{Definition}
\newtheorem{remark}{Remark}

\newenvironment{pf}
{\medskip\noindent {\it Proof --- \ }}
{\hfill\nobreak $\Box$ \par\bigbreak}

\newcommand{\RR}{{\mathcal R}}

\newcommand{\Hom}{\text{Hom}}

\newcommand{\C}{{ \mathbb C  }}
\newcommand{\R}{{ \mathbb R  }}
\newcommand{\Q}{{ \mathbb Q } }
\newcommand{\Z}{{ \mathbb Z  }}

\newcommand{\HH}{{\mathcal H}}

\newcommand{\Gal}{{\mathrm{Gal}\,}}

\newcommand{\A}{\mathbb A}

\newcommand{\anneau}{{ \mathcal O}}

\newcommand{\U}{{\text{U}}}

\newcommand{\Gl}{{\text {GL}}}

\newcommand{\GL}{{\text {GL}}}

\newcommand{\Hecke}{{\mathcal{H}}}

\newcommand{\G}{{\mathfrak g}}

\newcommand{\ch}{{\text{ch}}}

\newcommand{\ses}{{\text{ss}}}

\newcommand{\tr}{{\text{tr\,}}}

\newcommand{\rig}{{\text{rig}}}

\newcommand{\cal}{\mathcal}

\newcommand{\W}{W_{\Q_p}}

\renewcommand{\Gal}{{\rm Gal}}

\renewcommand{\AA}{\mathcal{A}}

\renewcommand{\G}{G}

\newcommand{\crys}{\text{crys}}

\begin{document}

\baselineskip 14.8pt

\title{New examples of $p$-adically rigid automorphic forms}

\thanks{This note arose from a question that Laurent Clozel asked me 
in 1997, when I was a young graduate student under his supervision. 
The question was in substance: 
``the representation $\prod_{v \text{place of $\Q$}} \pi^n(\chi)$
is cuspidal for $\U(2,1)$ if and only if $L(s,\chi)$ vanishes at $s=1/2$ at an even order. Is it possible to infer some 
arithmetically significant informations from this fact?''.
Let this be the occasion to thank Laurent Clozel for this question, and also
for all his help during my  PhD and thereafter. 
I also thank Ga\"etan Chenevier, with which I had innumerable conversations
on subjects related to this article, for his patience, encouragement, and help.
During the writing of this article, I was supported by the NSF grants 
DMS 05-01023 and DMS 08-01205}

\author[J.~Bella\"iche]{Jo\"el Bella\"iche}
\email{jbellaic@brandeis.edu}
\address{Jo\"el Bella\"iche\\Brandeis University\\
415 South Street\\Waltham, MA 02454-9110\\U.S.A}

\maketitle

\bibliographystyle{style} 


\section{Introduction}

We denote by $\U(2,1)$  the quasi-split unitary group 
in three variables attached to a quadratic imaginary field,
and we fix a prime $p$ at which $\U(2,1)$ splits.
  
 \par \bigskip 
The main aim of this note is to show (Theorem~\ref{rig}) that automorphic 
representations for $\U(2,1)$ that are cohomological in degree 1, together with a suitable choice of a refinement\footnote{Some authors, like \cite{SU}, call
the choice of a refinement a $p$-{\it stabilization}} at $p$
  are {\it $p$-adically rigid}, i.e. cannot be
put in any non-trivial positive-dimensional $p$-adic family. 

This result is  a little bit surprising, at least for me, since, as we shall see, the 
set of weights of those automorphic representations is not $p$-adically discrete, and is in fact dense
in a 1-dimensional\footnote{In all this introduction, we are ruling out central (that is, anticyclotomic) deformations. 
Otherwise, the dimension would be 2.} subspace of the weight space. Furthermore, the same automorphic representations with another choice of refinement
can easily be put in a non-trivial 1-dimensional $p$-adic family 
(see Remark~\ref{semiordinary}).
Hence the eigenvariety for $\U(2,1)$ and cohomological degree $1$, if it exists
in any reasonable sense (e.g \cite{E}), is not equidimensional. 
Even more, different refinements (or $p$-stabilizations)  
of a same automorphic form may lie in components of different 
dimensions of the eigenvariety.

A rigidity result for some examples of ordinary, 
non essentially self-dual cohomological automorphic forms
of $\Gl_3$ has already been obtained by Ash, Pollack, and Stevens~\cite{APS}.
They also conjecture that all non essentially self-dual automorphic forms for 
$\Gl_n$, $n \leq 3$ are rigid. Their is no doubt that their  result is 
much deeper, both in significance and in difficulty, 
but the point of this paper is to show that different kind of
pathologies may occur in the case of unitary groups.
 
Non self-dual cohomological automorphic forms for $\Gl_n$ are very mysterious.
If $n \geq 3$, there is no Shimura varieties for $\Gl_n$ or any of its inner form, and no one knows how to attach a 
Galois representation to a non self-dual cohomological automorphic forms for $\Gl_n$. Similarly, $\Gl_n(\R)$ 
has no discrete series for $n \geq 3$ and no one ever constructed 
a non-trivial, positive-dimensional families for those groups except by using functoriality from some other groups. The rigidity result of \cite{APS}
for some automorphic forms for $\Gl_3$ fits well in this picture.

The case of $\U(2,1)$ is very different. Hida constructed 
nice, equidimensional (of the expected dimension) families in the ordinary 
case for this group and many other a long time ago,
and we also have such families for the definite inner form $\U(3)$ of $\U(2,1)$ (and higher rank analogs) without the
ordinarity assumption, by works of Chenevier~\cite{C} and Emerton~\cite{E}. Urban and Skinner
claim in \cite{SU} that they are able to construct such families 
in the semi-ordinary case\footnote{
Unfortunately, they do not give their proof.  Indeed, almost all non-trivial results
of \cite{SU} can be put,
according to the justification that is given for them, in the following tripartite classification:
\begin{itemize} 
\item[(a)] Results for which a (sound) reference is given, but of which no proof can be found in the given reference (e.g. loc. cit., Prop. 4.2.2(i)).
\item[(b)] Results for which "a proof will appear elsewhere", and for which no proof
has appeared anywhere (not even in a preprint, twenty months after the {\it publication} of \cite{SU}).
(See e.g. Theorem 3.2.1 loc. cit.) 
\item[(c)] Results for which a reference is given, that does not exist. Indeed, many results
are said to be proved in [U06], which refers to "Urban, E., Eigenvarieties for reductive groups. Preprint, 2006." This preprint does not exist, either in paper or electronic form, under this name or any other,
and for obvious reasons will never do.
\end{itemize}}, at least for representations that are holomorphic discrete series 
at infinity. This note shows that in the non-definite, non ordinary case, 
pathologies may occur, and that beyond a general theorem on the existence
of eigenvarieties (like, e.g. \cite{E}), a formula computing the dimension of their components
containing a refined automorphic representation would be 
very useful, and that this formula should depend not only of
the given representation, but also of the chosen refinement.
 
 \par \bigskip
   
 The discussion above was essentially concerned  by
 families of automorphic forms that are cohomological in a fixed degree -- as in the works
 of Stevens, Chenevier,  Emerton, and others. Our second result 
(Theorem~\ref{rig2}) deals with
  families where forms are cohomological, but not necessarily of the same degree. Our result is that through some
 very specific example of non-tempered representation of $\U(2,1)$, such more general families do not exist (more precisely, have dimension $0$.)

\section{Preliminaries}

\subsection{Unitary groups}

Let $E$ be a quadratic imaginary field, and let $G=\U(2,1)$ be the unique
quasi-split unitary group in three variables over $\Q$ that splits in $E$.

Let $p$ be a prime number that, we shall assume, splits in $E$. We shall also assume that $p \neq 13$. I don't think this is really useful, but who knows?

Let $K^p=\prod_{l \neq p} K_l$ be a compact open subgroup of $\G(\A_f^p)$ that we will call the {\it tame level}. We denote by $S$ the finite set of primes $l$
such that either $l=p$, or $l \neq p$ and $K_l$ is not maximal hyperspecial.
For $l \not \in S$, we denote by $\HH_l$ the Hecke algebra (over $\Z$) of $\G(\Q_l)$ with respect to $K_l$. 
We denote by $K_p$ a maximal compact hyperspecial subgroup of $G(\Q_p)$, by $I_p$ an Iwahori subgroup of 
$K_p$ and by $\AA_p$ the Atkin-Lehner sub-algebra of the Hecke algebra
of $G(\Q_p)$ with respect to $I_p$ -- See \cite[\S6.4.4]{BCbook}. 
We set $\HH = \left( \bigotimes_{l \not \in S} \HH_l \right) \otimes 
\AA_p$. This is a commutative algebra over $\Z$.

\subsection{Automorphic refinements}
We refer to \cite[\S5.4]{BCbook} for a complete discussion of  
the notion of {\it refinement} $\RR$ of 
an unramified  representation $\pi_p$ of $G(\Q_p)$. Let us simply define it here as a character of $\AA_p$ 
appearing in $\pi_p^{I_p}$.

\subsection{Hecke characters} \label{unitary}

Set $K=K_p \times K^p$.
A {\it refined representation} $(\pi,\RR)$ of level $K$ will mean a discrete
automorphic representation $\pi$ for $\U(2,1)$, such that $\pi^K \neq 0$,
together with a refinement $\RR$ of $\pi_p$. Such a refined representation
defines a character $\psi_{\pi,\RR} = \psi_\pi \otimes \psi_\RR: \HH 
\rightarrow \C$. Here, $\psi_\pi : \left( 
\bigotimes_{l \not \in S} \HH_l \right) \rightarrow \C$ is the character
giving the action of $\left( 
\bigotimes_{l \not \in S} \HH_l \right)$ on the one-dimensional space
$\prod_{l \not \in S} \pi_l^{K_l}$, and $\psi_\RR$ is the character of $\AA_p$ defined by the refinement $\RR$.

If $\pi_\infty$ is 
cohomological (for some coefficient system), the image of $\psi_{\pi,\RR}$
lies in fact in a number field. We fix once and for all an embedding $\bar \Q \subset \bar \C$, and also an embedding $\bar \Q \subset \bar \Q_p$. Hence we
can see $\psi_{\pi,\RR}$ as a character $\Hecke \rightarrow \bar \Q_p$ when 
$\pi_\infty$ is cohomological. Conversely, we shall say that a character 
$\psi_0 : \HH \rightarrow \bar \Q_p$ is of level $K$ and cohomological of 
degree $n$ if it is of the form $\psi_{\pi,\RR}$ for $(\pi,\RR)$ a refined 
representation of level $K$ such that $\pi_\infty$ is cohomological 
in degree $n$, and $\pi^K \neq 0$.

\subsection{Galois representations}

\label{galoisrep}

To a discrete automorphic representation $\pi$ such that $\pi_\infty$ is 
cohomological of some degree is attached by the work of Blasius and Rogawski
(\cite[Thm 1.9.1]{BR}) and the whole book \cite{ZFPMS} a three dimensional representation (see  
\cite[\S3.2.2]{BC} for this formulation of the results of Blasius and Rogawski) 
$\rho_\pi : \Gamma_E \rightarrow \Gl_3(\bar \Q_p)$ that depends only of the 
Hecke-character  $\psi_{\pi}$ (Here and below, 
$\Gamma_F$ denotes ``the'' absolute Galois group of a field $F$)

The representation $\rho_\pi$ satisfies the symmetry condition $\rho_\pi = \rho_\pi^\bot$, where for $\rho$ a 
representation of $\Gamma_E$, $\rho^\bot(g) = {}^t \rho(cgc)^{-1}$. Here,
$c$ is an element of $\Gamma_\Q-\Gamma_E$ of order $2$.
 
We fix once and for all a place $v$ of $E$ above $p$. The decomposition group $D_v$ at $v$ in $\Gamma_E$ is isomorphic to $\Gal(\bar \Q_p/\Q_p)$. The choice of $v$ determine an isomorphism $G(\Q_p) \simeq \Gl_3(\Q_p)$, well-defined up to conjugation.

If $\pi_p$ is unramified, then $\rho_{\pi}|D_v$ is crystalline. It thus
has three integral Hodge-Tate weights, that we will denote by
$k_0 \leq k_1 \leq k_2$ and three Frobenius eigenvalues $\phi_0,\phi_1,\phi_2$.
 The weak admissibility condition
of crystalline representation theory implies the equality 
$v_p(\phi_0)+v_p(\phi_1)+v_p(\phi_2) = k_0+k_1+k_2$, where $v_p$ is the $p$-adic valuation.

We will refer to the $k_i$'s and the $\phi_i$'s simply as the weights and eigenvalues of $\rho_\pi$, 
without mentioning the choice of the place $v$ above $p$. 

\subsection{Refinement on the Galois side}

For our purposes, a {\it refinement} of $\rho_\pi$ is simply an ordering of the eigenvalues
$\phi_0,\phi_1,\phi_2$. 
In a generic situation, there is thus 6 refinements, and less
if two eigenvalues are equal. The important fact for us (see \cite[\S6]{BC}) 
is that each automorphic refinement at $p$ of $\pi$ defines canonically a Galois refinement of 
$\rho_\pi$. Refinements of $\rho_\pi$ that can be obtained this way  
are called {\it accessible}.

\subsection{The weight space}

We define a weight space $\W := \Hom_{\rig}(\Z_p^\ast,G_m)$ of dimension $2$,
and embed $\Z^2$ as a Zariski-dense subspace in $\W(\Q_p)$ in a natural way.

\subsection{Representations that are cohomological in degree $1$}

\label{rog}

Let $\chi_0$ be a Grossencharacter of $E$ such that 
for all $z \in \AA_E$, $\chi_0(z \bar z)=1$ and such that for 
$z \in E\otimes R = \C$, $(\chi_0)_\infty(z)=z^k/(z \bar z)^{k/2}$, where $k$ is an odd
positive integer. (The oddness of $k$ is equivalent to the property that
$\chi_0$ does {\it not} come from a Hecke-character of $\U(1)$.) 

The complete $L$-function $L(\chi_0,s)$ satisfies a functional equation (see \cite[3.6.8 and 3.6.1]{Ta})
$$L(\chi_0,1-s) = \epsilon(\chi_0,s) L(\chi_0,s),\ \ \epsilon(\chi_0,1/2)=\pm1 .$$

For every place $v$ of $\Q$, we can define an irreducible admissible 
representation $\pi^n(\chi_0)_v$ of $G(\Q_v)$ which is non tempered 
(and unramified if $v$ is finite, unramified in $E$, and $(\chi_0)_v$ is unramified).
The representations $\pi^n(\chi_0)_v$ belongs to only one local $A$-packet), 
which is a singleton when $v$ is split in $E$, and a pair $\{\pi^n(\chi_0)_v,\pi^s(\chi_0)_v\}$
when $v$ is ramified or inert in $E$ 
(including the non-archimedean place of $\Q$).
A representation of the global corresponding $A$-packet $\Pi(\chi_0)$ is
of the form $\pi = \otimes_v \pi_v$, where $\pi_v = \pi^n(\chi_0)_v$ 
for almost all $v$, and $\pi_v=\pi^s(\chi_0)_v$ otherwise. Let $n(\pi)$ be the number of places $v$ of $\Q$ 
such that $\pi_v=\pi^s(\chi_0)_v$.

\begin{fact}[Rogawski]  \label{fact1}
\begin{itemize} 
\item[(i)] Assume that $\chi_0$, $\pi$ are as above and that
 $(-1)^{n(\pi)} = \epsilon(\chi_0,1/2)$.
Then $\pi$ is discrete, and cohomological in degree $1$ (and $3$).
Moreover $\pi$ is cuspidal, except in the case where $n(\pi)=0$ and $L(1/2,\chi_0) \neq 0$.
\item[(ii)] Conversely, any discrete automorphic representation of $\G$
that is cohomological in degree $1$ is of the form described in (i).
\end{itemize}
\end{fact}
Reference for this fact are: \cite{Rmult} for the definition of $\pi^n$, $\pi^s$ and the existence assertion;
\cite{Rbook} for the cuspidality assertion, and \cite[4.4]{R} for the cohomological assertions of (i), and (ii). (For (ii), see also \cite[Prop 3]{MR} together with the classification recalled in \cite[\S3]{BC}).

\par \bigskip

The character $\chi(z) = \chi_0(z) (z\bar z)^{1/2}$ is algebraic. 
We still denote by $\chi$ its $p$-adic realization as a character $\Gamma_E \rightarrow \bar \Q_p$.
We have $\chi^\bot=\chi(1)$.

\begin{fact}[Rogawski]\label{fact2}
If $\pi$ is as in the above fact, then $\rho_\pi$ is the sum of three 
characters $\chi \oplus 1 \oplus \chi^\bot$, and 
the weights $k_0 \leq k_1 \leq k_2$ of $\rho_\pi$
satisfy either $k_1 - k_0 = 1$ or $k_2 - k_1 = 1$
\end{fact}
Indeed, the description of $\rho_\pi$ follows directly form the description of $\pi^n$ in \cite{Rmult}
(see also \cite[\S3.2.3]{BC}), and it follows that (if $v$ is well chosen) the weights of $\rho_\pi$
are in some order $0, (k-1)/2, (k+1)/2$, which implies that $k_2-k_1=1$ if $k>0$, and $k_1-k_0=1$ if $k<0$.

 \begin{remark} It is easy to see that if the level $K^p$ is small enough,
then for all triples of integers 
$(k_0,k_1,k_2)$ that satisfies either $k_1-k_0=1$, or $k_2-k_0=1$,
there are exists representations $\pi$ as in Fact ~\ref{fact1} with weights 
$(k_0,k_1,k_2)$.
The Zariski closure of those weights in the three dimensional weight space obviously has dimension $2$,
or if we fix one of the weight, for example the first one $k_0$ to $0$, has dimension $1$.
So it is natural to expect that there exists one-dimensional non-trivial
$p$-adic families interpolating forms that are cohomological
in degree $1$. Theorem~\ref{rig} shows that actually it is not 
always the case.
\end{remark}

Let $\pi$ be as in Fact~\ref{fact1}. Assume in addition that 
$\pi_p$ is unramified. By Fact~\ref{fact2}, $\rho_\pi$ is the sum of three characters, so
the weak admissibility conditions applied to each of them gives that the
set of weights $\{k_0,k_1,k_2\}$ is equal to the set of slopes 
$\{s_1,s_2,s_3\}$ of $\rho_\pi$. Let us assume that the weights are distinct 
(that is, that $|k| > 1$) to simplify the discussion,
 Thus, a refinement of $\rho_\pi$ is an ordering of the weight $\{k_0,k_1,k_2\}$. 

\begin{fact} \label{fact3} 
If $\pi_p$ is unramified, then $\pi_p^{I_p} = \pi^n(\chi_p)^{I_p}$ has dimension $3$, and $\pi$ has three possible refinements. Assume
to fix ideas that $k_2-k_1=1$ (the other case being symmetric). 
The corresponding accessible 
refinements of $\rho_\pi$ are (a) $(k_0,k_2,k_1)$; 
(b) $(k_2,k_0,k_1)$, and (c) $(k_2,k_1,k_0)$
\end{fact}
This follows from \cite[Remark 5.2.4]{BC}.

\begin{remark} 
\begin{itemize}
\item  The {\it ordinary} refinement would be $(k_0,k_1,k_2)$. It is not 
{\it accessible} in the present situation. In other words, we are not in a case where Hida's work allow us to construct 
$p$-adic families through $\pi$.

\item The refinement $(k_2,k_0,k_1)$ is {\it anti-ordinary} in the sense of 
\cite[\S2.4]{BCbook}. It is the refinement used in the paper \cite{BC}.

\item The refinement (a) is, I guess, what Urban and Skinner called 
{|it semi-ordinary}. It is analog to the one they
use in their paper \cite{SU2}.

\item I do not see anything to say on refinement (c).
\end{itemize}
\end{remark}

\section{Rigidity of cohomological in degree 1 forms.}

\subsection{Definition of a $p$-adic family of forms that are cohomological in degree $n$}

For our purpose, we will adopt the following, rather weak, 
notion of $p$-adic family. Let $K$ be as in~\S\ref{unitary}

\begin{definition} \label{fam} 
Let $n$ be an integer. By a 
{\it $p$-adic family interpolating discrete automorphic representations of level $K$ and cohomological of degree $n$} for $\U(2,1)$,
or for short {\it a $p$-adic family}, we mean the data of

\begin{itemize}
\item A rigid analytic space $X$ over $\Q_p$ that is reduced and separated.
\item An analytic map $\kappa = (\kappa_1,\kappa_2) : X \rightarrow \W$ 
over $\Q_p$.
\item A ring homomorphism $\psi : \HH \rightarrow \anneau(X)^\rig$.
\item A subset $Z \in X(\bar \Q_p)$ that is Zariski-dense and accumulates 
at every of its points.
\end{itemize}

such that
\begin{itemize}
\item[(a)] For every $z \in Z$, the character $\psi_z : 
\Hecke \rightarrow \anneau(X)^\rig \stackrel{\text{eval. at $x$}}{\rightarrow} 
\bar \Q_p$ is of the form $\psi_{\pi,\RR}$ where $\pi$ is of level $K$ and cohomological of degree $n$. 
\item[(b)] if $z \in Z$, and $\psi_z = \psi_{\pi,\RR}$, then the Hodge-Tate
weights of $\rho_\pi$ are the natural integers $0,\kappa_1(z),\kappa_2(z)$, 
and we have $0 \leq \kappa_1(z) \leq \kappa_2(z)$.
\item[(c)] The map $Z \rightarrow \Hom(\HH, \bar \Q_p),\ z \mapsto \psi_z$ 
is injective.
\end{itemize}
\end{definition}
It is expected that eigenvarieties, if they exist, satisfy 
those properties, and many more: see 
\cite[Definition 7.2.5 and theorem 7.3.1]{BC} or \cite[Theorem 0.7]{E}.

\begin{remark}
The condition (c) is a condition of non-triviality, though very weak. It ensures that the family is not the constant family, for example. Formulating sensible non-triviality conditions is not so easy. This difficulty is related 
is related to the one one have to  interpret eigenvarieties as  moduli spaces.

\end{remark}

\subsection{Rigidity}

We shall say that a given refined automorphic representation {\it belongs} to
the family $X$ if there is a $z \in Z$ such that $\psi_z = \psi_{\pi,\RR}$.

\begin{theorem} \label{rig} Let $X$ be an {\it irreducible} $p$-adic family as in 
definition~\ref{fam} of cohomological degree $n=1$. Suppose that a refined 
discrete representations $\pi$ with distinct weights and 
with a refinement $\RR$ of type (b) or (c) of Fact~\ref{fact3} 
belongs to $X$. Then $X$ is a point.
\end{theorem}

\begin{prop}  \label{gal} Let $X$ be a family as in the definition~\ref{fam}. 
There exists a pseudo-character $T : \Gamma_E \rightarrow \anneau(X)^\rig$ and three locally constant functions $F_0,F_1,F_2$ on $X(\bar \Q_p)$ such that 
for every $z \in Z$, such that if $\psi_z = \psi_{\pi,\RR}$,
and $\rho_z$ is the semi-simple representation of trace $T_z := \text{eval}_z 
\circ T$, then $\rho_z \simeq \rho_\pi$, 
and the eigenvalues  or $\rho_\pi$
are $F_0(z), F_1(z)p^{\kappa_1(z)}, F_2(z)p^{\kappa_2(z)}$, given 
in the order defined by the refinement $\RR$.
\end{prop}
\begin{pf}
See \cite[\S7]{BC}. The functions $F_0,F_1,F_2$ are defined as the 
images by $\psi$ of three 
suitable elements in the Atkin-Lehner algebra $\AA_p$. 
\end{pf}

Let us now assume, by contradiction, that $\dim(X) > 0$. By hypothesis,
every $z$ in $Z$ correspond to a $\pi$ that is cohomological in degree $1$.
so the weights of $\pi$, that is by (b) of the definition of a family
$0 \leq \kappa_1(z) \leq \kappa_2(z)$ satisfy, by Fact~\ref{fact2}, 
either $\kappa_1(z)=1$ or $\kappa_2(z)-\kappa_1(z)=1$. Since 
$X$ is irreducible, and $Z$ Zariski-dense, one of those two equalities has 
to hold for all $z \in Z$ (and even over $X$). By symmetry, we may and do 
assume that we have for all $z \in Z$, $$\kappa_2(z) - \kappa_1(z)=1.$$

Now, by hypothesis, there exists a $z_0$ in $Z$ such that the refinement
is not the refinement (a). That is to say, the refinement at $z_0$ is either 
$(\kappa_2(z_0), \kappa_1(z_0),0)$ or 
$(\kappa_2(z_0),0,\kappa_1(z_0))$. Using Prop~\ref{gal} we 
see that in both cases $v_p(F_0(z_0))=\kappa_2(z_0)$.

Let $U$ be an affinoid neighborhood of $z_0$ on which the function $v_p(F_2)$ 
is constant. By hypothesis $Z \cap U$ is infinite.
Since there is only a finite number of automorphic representations with a 
fixed level and fixed weight, there exists $z \in Z \cap U$ such that 
$\kappa(z) \neq \kappa(z_0)$. This implies 
$\kappa_2(z) \neq \kappa_2(z_0)$, since $\kappa_1 = \kappa_2 -1$ everywhere.
But by Prop~\ref{gal} applied to $z$, we have $\kappa_2(z) = v_p(F_2(z)) = v_p(F_2(z_0)) = \kappa_2(z_0)$, a contradiction. This proves Theorem~\ref{rig}.

\begin{remark} \label{semiordinary} Starting with a $1$-dimensional
deformation $\chi_x$ (for $x$ a parameter in a rigid analytic space $X$)
of the Hecke character $\chi$, it is very easy to construct
a one dimensional family as above of representation that are cohomological
of degree $1$, refined with a refinement of type (a). This may be done
``by hand'', defining explicitly the function $\psi(h)$, $h \in \HH$, with the guidance that $\rho_x = 1 \oplus \chi_x \oplus \chi_{x}^\bot$. We leave the
details to the reader.
\end{remark} 

\begin{remark} For higher rank unitary groups over $\Q$ 
(of signature $(n-1,1)$ at infinity), it should be possible to exhibit similar
examples
of rigidity in cohomological dimension $1$, using instead of Fact~\ref{fact2} 
the main theorem of \cite{MR}, but we have not written down the details.
\end{remark}

\section{Another rigidity result}

In this section, we give another rigidity result, which holds for a 
much smaller class of representations, but for a more general notion of family, 
where we do not assume the classical representations 
 to have a fixed cohomological degree.

\subsection{A special non-tempered refined representation}

Let $\chi_0$ be Hecke character of $E$ as in \S\ref{rog}, and assume that $L(\chi_0,1/2) \neq 0$. 
In particular, $\epsilon(\ch_0i,1/2)=1$, and by Fact~\ref{fact1} 
there exists a discrete
automorphic representation $\pi^n(\chi_0)=\otimes_{v} \pi^n(\chi_0)_v$ 
for $G=\U(2,1)$. Let $\RR$ be the refinement (b) of $\pi^n(\chi_0)_p$ (see Fact~\ref{fact3}).

\subsection{Definition of a family}

\label{special}

We fix a subgroup $K$ of $\U(2,1)(A_\Q)$ and a set $S$ of bad primes
as in~\ref{unitary}. We set $S^{p} = S-\{p\}$.

Since the argument we will give is arithmetic. We will need to control
the ramification at every places so we shall had a condition at primes 
of $S^{p}$. This can be done by imposing some types (in the sense 
of Bushnell and Kutzko) at those primes as in \cite{BC} or more 
conveniently, by using the notion of NMRPSS as in \cite[\S6.6]{BCbook}.
We choose the second way. Note that the representation $\pi^n(\chi_0)$
is an NMRPSS at every place.

\begin{definition} \label{fam2} 
By a 
{\it strong $p$-adic family interpolating discrete cohomological automorphic 
representations of level $K$} for $\U(2,1)$,
or for short {\it a strong $p$-adic family}, we mean the data of

\begin{itemize}
\item A rigid analytic space $X$ over $\Q_p$ that is reduced and separated.
\item An analytic map $\kappa = (\kappa_1,\kappa_2) : X \rightarrow \W$ 
over $\Q_p$.
\item A ring homomorphism $\psi : \HH \rightarrow \anneau(X)^\rig$.
\item A subset $Z \in X(\bar \Q_p)$ that is Zariski-dense and accumulates 
at every of its points.
\end{itemize}

such that
\begin{itemize}
\item[(a)] For every $z \in Z$, the character $\psi_z : 
\Hecke \rightarrow \anneau(X)^\rig \stackrel{\text{eval. at $z$}}{\rightarrow} 
\bar \Q_p$ is of the form $\psi_{\pi,\RR}$, where $\pi$ is of level $K$, cohomological of degree $n$, and $\pi_l$ is an NMRPSS 
(see \cite[Def 6.6.5]{BCbook}) for all $l \in S^p$.
\item[(b)] if $z \in Z$, and $\psi_z = \psi_{\pi,\RR}$, then the Hodge-Tate
weights of $\rho_\pi$ are the natural integers $0,\kappa_1(z),\kappa_2(z)$,
and we have $0 \leq \kappa_1(z) \leq \kappa_2(z)$.
\item[(c)] The map $Z \rightarrow \Hom(\HH, \bar \Q_p),\ z \mapsto \psi_z$ 
is injective.
\end{itemize}
\end{definition}

\subsection{The Kisin property}

The proof of Prop~\ref{gal} works without any change for a family
$X$ as above. We keep its notations. For $z \in Z$, we denote by $\anneau_z$
the rigid analytic local field of $X$ at $z$.

\begin{definition} A {\it realization} of $T$ over $\anneau_z$ 
is any torsion-free finite module $M$ over $\anneau_{z}$ with
a continuous action of $\Gamma_E$, such that for all $g \in \Gamma_E$,  $\tr(g | M \otimes \text{Frac}(\anneau_z))$ is in $\anneau_z$ and is equal 
to the germ of $T(g)$ at $z$.
\end{definition}

\begin{definition} We shall say that the family $X$ satisfies the
 Kisin property
at $z \in Z$ if for any realization $M$ of $T$ over $\anneau_z$
we have, if $\dim D_\crys(M_z^{\ses})^{\phi=F_0(z)} \leq 1$, then
$\dim D_\crys(M_z)^{\phi=F_0(z)} = 1$.

We say that the dual Kisin property holds if the same implication 
holds with $T$ replaced by the dual character $T(g^{-1})$ and $F_0(z)$ by $p^{-\kappa_2(z)} F_2(z)$.

We say that the strong Kisin property holds if both the Kisin and the dual Kisin properties hold.
\end{definition}

It is expected that eigenvarieties, and most of the naturally constructed families satisfies this
properties. For the case of the eigencurve, this is essentially proved in \cite{kisin},
and for the eigenvarieties on unitary groups, and more generally families with a
sufficiently dense set of "classical" points, in \cite[\S3.4]{BCbook}.

\begin{remark}

The reason for which I state the Kisin property 
as a definition rather than
as a theorem is that I do not want to presume what are the minimal 
conditions on a family that implies it. 
As I have just said, this properties has been proved in \cite{BCbook}
under some assumptions on the set $Z$, and a similar property is stated under 
a different assertion on $Z$ in \cite[Prop. 4.2.2(i)]{SU}, without proof -- the reference given there  to \cite{kisin} must be a joke. Anyway, since we want to prove the 
most general possible rigidity result in its line of reasoning, 
it is better to assume directly the Kisin property we need rather than any property on the family 
that implies it. 

Another remark is that the formulation of the Kisin property given above,
for general torsion free (as opposite to free, or locally free) realizations,
may seem weird and complicated. Indeed, this formulation has been chosen 
merely for convenience. Let me explain the issue at stake. 
In general, it is not possible to choose in the neighborhood of a reducible
point $z_0$ a realization of the pseudo-character $T$ that is free 
(see the discussion in \cite[Chapter 1]{BCbook} and \cite[Remark 9]{B}).
However, as far as we only want a Kisin's property on the representation
at a closed point, we may in practice  perform some surgery on $X$ to
get a free module. One way to do so, simple though unestehtic, 
is to replace $X$  by the normalization on a suitable curve in $X$ through $z_0$.
(This is the way used in~\cite{BC}, and later in~\cite{SU2}). Another way
is to replace $X$ by a suitable blow-up (a normalization may not be enough,
as shown by the example developed in~\cite{B}, since the normalization of the family
defined in \cite[\S5]{B} still satisfies the properties (i) to (ix) loc. cit.), idea which is used in 
\cite[Chapter 3]{BCbook} (and before, though for a different purpose, in \cite{kisin}). Now, those surgical operations work for real-world 
families, but in an abstract setting like ours 
(or the one of \cite{SU}), one would need a definition of family which 
is stable by those operations. This would 
lead to unnecessarily restrictive and  unesthetic  assumptions in the 
definition of a family. So we have chosen another way, namely a
 formulation of the Kisin property for general torsion-free modules, which is not 
much harder to prove in practice. 

Let us also note that the way Urban and Skinner in 
\cite{SU} deal with this issue is doubly mistaken - though it can be easily 
rectified, along the way explained above. Indeed they only formulate 
their version of the Kisin property (\cite[Prop. 4.2.2(i)]{SU}) for free modules. 
But when they apply it,
page 496, they unduly assume that their module ${\cal L}$, 
over a two-dimensional normal basis, is locally free\footnote{I am not saying
that in their specific situation ${\cal L}$ is not locally free. I am simply saying that this fact does not follow from what they say, and that they have to prove it.} Moreover, since they apply (page 496) Prop. 4.2.2(i) to a 
sufficiently small affinoid, they implicitly
assume  that the restriction of a finite slope family in their sense 
(loc. cit. 4.2) to an 
affinoid subset is still a finite slope family -- which is false. 
\end{remark}

\subsection{Rigidity}

\begin{theorem} 
\label{rig2}
Let $X$ be an {\it irreducible} family as in Definition~\ref{fam2},
and assume that $(\pi^n(\chi_0),\RR)$ (defined in \ref{special}) belongs to $X$, and that $X$ satisfies the strong Kisin property at the corresponding point $z_0$. Then $X$ is a point.
\end{theorem}
\begin{pf} 
Let $Z_1$ be the set of points of $Z$ which are cohomological of degree $1$.
If $Z_1$ is Zariski-dense, then replacing $Z$ by $Z_1$ we see that $X$ is a family in the sense of Definition~\ref{fam}, and $X$ is a point by 
Theorem~\ref{rig}. Therefore we
can assume that $Z_1$ is not Zariski dense.

Let $Z_0$ be the set of points such that $\psi_z = \psi_{\pi,\RR}$ for $\pi$ a one dimensional representation. We have at those points 
$\rho_z = \mu_1 \oplus \mu_2 \oplus \mu_3$, with the three characters $\mu_i$ satisfying $\mu_i = \mu_i^\bot$.
If $Z_0$ was Zariski-dense, the same property would hold for any $\rho_z$, $z \in Z$.
But this is false for $\rho_{z_0} \simeq \chi \oplus 1 \oplus \chi^\bot$, since $\chi \neq \chi^\bot =\chi(1)$.
Therefore, $Z_0$ is not Zariski-dense.

Hence we see that there is a Zariski-dense set of $z$ corresponding to points 
that are neither cohomological of degree $1$, nor $1$-dimensional. By 
Rogawski's classification (see \cite[\S 3.2.3]{BC}), those points have a $\rho_z$ that are 
either irreducible, or the sum of a character and a two-dimensional 
representation. Thus the 
pseudo-character $T$  either is irreducible, or is the sum of a 
pseudo-character of dimension $2$ and a character.

If $T$ is the sum of a pseudo-character $T'$ of dimension $2$ and of 
a character, then $T'_{z_0}$ is stable by the operation $\bot$, so we have 
$T'_{z_0} = \chi \oplus \chi^\bot$. From this, we can construct 
a non trivial extension $U$ of $\chi$ by $\chi^\bot$ 
(as a sub-quotient of some realization $M$ of $T$ over $\anneau_{z_0}$, with 
$D_\crys(M^\ses_{z_0})^{\phi=F_0(z_0)}$ of dimension $1$ -- see 
\cite[Theorem 1.5.6]{BCbook}) that satisfies $U \simeq U^\bot$. 
From the Kisin properties at $z_0$, $U$ is crystalline at $v$. Since $U \simeq U^\bot$ it is also 
crystalline at $\bar v$. The twist $U \chi^{-1}$ is an extension of $1$ by $\bar \Q_p(1)$
which is crystalline at $v$ and $\bar v$, and non ramified everywhere else 
by the level $K$ and NRMPSS assumption (see \cite[Prop 8.2.10]{BCbook}).
As is well known, there is no such extension in the category of $\Gamma_E$-representation, $E$ a quadratic imaginary field (see e.g. 
\cite[Prop. 5.2.2]{BCbook}). This is a contradiction.

If $T$ is irreducible, then let us simply say that we can 
proceed exactly as in  as in 
\cite[\S9]{BC} (using our strong Kisin hypothesis instead of the 
local version of Kisin as in \cite[Chapter IX]{BCbook}) 
to construct an element in $H^1_f(E,\chi)$.
But since $L(\chi_0,1/2) \neq 0$, Bloch-Kato conjecture and a theorem of 
Rubin (\cite{rubinmc})  agrees that such an extension does not exist. 
 \end{pf}

\begin{remark} It would be interesting to confront the above rigidity 
result with 
\cite[Theorem 3.2.1]{SU}, which states the existence of  
two-dimensional families for unitary groups in some cases, and its proof,
at least when the authors provide a precise general statement and a complete argument  proving it. As it stands, the sketch of proof of \cite[Theorem 3.2.1]{SU} is so sketchy 
that it applies smoothly to our representation $\pi^n(\chi_0)$, 
leading to a contradiction with the rigidity theorem~\ref{rig2}.

\end{remark}

\par \bigskip

\end{document}